\documentclass{article}
\usepackage{fullpage}
\usepackage{graphics,psfrag,epsfig}
\usepackage{amsfonts,latexsym,eucal,amsmath,amsthm,amssymb}

\begin{document}

\newcommand{\rum}{\rule{0.5pt}{0pt}}
\newcommand{\rub}{\rule{1pt}{0pt}}
\newcommand{\rim}{\rule{0.3pt}{0pt}}
\newcommand{\numtimes}{\mbox{\raisebox{1.5pt}{${\scriptscriptstyle \times}$}}}
\newcommand{\optprog}[2]
{%
  \noindent\mbox{}\\[0cm]
  \noindent\fbox{%
  \begin{minipage}{0.955\linewidth}
    \mbox{}\\[-0.5cm]
    #1\\[#2]
  \end{minipage}
  }
  \noindent\mbox{}\\[-0.2cm]
}

\renewcommand{\refname}{References}

\twocolumn[%
\begin{center}
{\Large\bf Is it possible to compute the M\"{o}bius function without factoring? \rule{0pt}{13pt}}\par
\bigskip
Craig Alan Feinstein \\ {\small\it 2712 Willow Glen Drive, Baltimore, Maryland 21209\rule{0pt}{13pt}}\\
\raisebox{-1pt}{\footnotesize E-mail: cafeinst@msn.com, BS"D}\par
\bigskip
{\small\parbox{11cm}{%
\bigskip \noindent \textbf{Abstract:} It has been well known since Fermat's Little Theorem was first
published that it is possible to determine that a number is composite without determining any of its nontrivial factors. 
It is natural to ask whether it is also possible to compute the M\"{o}bius function of a composite number without determining any 
of its nontrivial factors. In this note, we argue that this is impossible.

\bigskip \noindent \textbf{Disclaimer:} This article was authored by Craig
Alan Feinstein in his private capacity. No official support or endorsement by the U.S. Government is
intended or should be inferred.\rule[0pt]{0pt}{0pt}}}
\bigskip
\end{center}]{%

It has been well known since Fermat's Little Theorem was first
published that it is possible to determine that a number is composite without determining any of its nontrivial factors \cite{b:FLT}. It 
is natural to ask whether it is also possible to compute the M\"{o}bius function of a composite number without determining any 
of its nontrivial factors; the M\"{o}bius function $\mu(n)$ is defined as $\mu(1)=1$, $\mu(n)=(-1)^k$ if $n=p_1 \cdots p_k$ 
where $p_j$ is a distinct prime for each $j=1,\dots,k$, and $\mu(n)=0$ if $n$ is divisible by a perfect square \cite{b:Mo}.
In this note, we argue that this is impossible:

Suppose that $n$ is square-free. Then $\mu(n)$ is defined as $(-1)^k$, where $k$ is the number of primes in the prime 
factorization of $n$. Then computing this expression requires computing $k \pmod 2$, and computing $k \pmod 2$ requires 
counting the prime factors of $n$. (Note that one does not have to count as $1,2,3,4,5,\dots$; one can count as 
$1,2,1,2,1,\dots$.) Hence, when $n$ is square-free, the only way to compute $\mu(n)$ is to count the prime 
factors of $n$, which requires determining its prime factors.

Now, suppose that $n$ is divisible by a perfect square. First, note that the only way to determine that a number 
is a perfect square is to find its square root. Then \textit{a fortiori}, the only way to determine that $n$ is 
divisible by a perfect square is to find the square root of a perfect square that divides $n$, which results in a nontrivial 
factor of $n$. Therefore, in order to compute $\mu(n)$ when $n$ is divisible by a perfect square, it is necessary to 
determine at least one of the nontrivial factors of $n$.

We can conclude that computing the M\"{o}bius function of a composite number is at least as difficult as finding a 
nontrivial factor of a composite number, unlike the problem of determining whether a number is composite.

}

\end{document}